\documentclass[12pt,a4paper]{article}

\usepackage{amscd}
\usepackage{amsfonts}
\usepackage{amsmath}
\usepackage{booktabs}
\usepackage{url}
\usepackage{hyperref}
\usepackage[small]{caption}

\newtheorem{theorem}{Theorem}

\newcommand{\F}{{\mathbb F}}

\begin{document}

\title{$q$-Steiner Systems Do Exist}
\author{{\sc Michael Braun}\footnote{University of Applied Sciences Darmstadt, \url{michael.braun@h-da.de}}\,\,  and {\sc Alfred Wassermann}\footnote{University of Bayreuth, \url{alfred.wassermann@uni-bayreuth.de}}}
\maketitle

\begin{abstract}
In this paper we give the first construction of a $q$-analog of a Steiner system. Using a computer search we found at least 26 $q$-Steiner Systems $S_2[2,3,13]$ admitting the normalizer of a singer cycle as a group of automorphisms.
\end{abstract}

\newpage
\section{Introduction}

A $t$-$(n,k,\lambda;q)$ \emph{design} over $\F_q$ is a set $\mathcal{B}$ of $k$-subspaces of the $n$-dimensional vector space $\F_q^n$ over the finite field $\F_q$ with $q$ elements such that each $t$-subspace of $\F_q^n$ is contained in exactly $\lambda$ elements of $\mathcal{B}$ \cite{BKL05,Ito98,MMY95,Suz90,Suz92,Tho87,Tho96}.

Designs over finite fields with parameters $t$-$(n,k,1;q)$ are called \emph{$q$-Steiner systems} and their parameters will be abbreviated by $S_q[t,k,n]$. In geometrical terms a $q$-Steiner system $S_q[t,k,n]$ is also called $(t,k)$-\emph{spread} in $\F_q^n$ \cite{HS11}. In case of $t=1$, $q$-Steiner systems are called \emph{spreads}. It is well known that a spread exists if and only $k$ divides $n$ \cite{LW01}.

Recently, in a seminal paper \cite{KK08} Koetter and Kschischang described an algebraic approach to the theory of error correcting subspace codes. A $q$-Steiner system $S_q[t,k,n]$ defines an optimal subspace code with minimum subspace distance of $2(k-t+1)$ \cite{KK08a}. 
Subspace codes defined by spreads are called \emph{spread codes} \cite{MGR08}.

Up to now, no $q$-Steiner systems $S_q[t,k,n]$ with parameter $t>1$ have been known. 
Metsch~\cite{Met99} even conjectured that such $q$-Steiner systems may not exist. 

The smallest admissible parameters of interest with $t>1$ for which $q$-Steiner systems $S_q[t,k,n]$ might be found are $q=2, t=2, k=3, n=7$.
Etzion and Schwartz \cite{SE02} discussed necessary and sufficient conditions for the existence of these $q$-Steiner systems. Moreover, Thomas \cite{Tho96} showed that certain $S_2[2,3,7]$ designs cannot exist. 

In \cite{EV11} Etzion und Vardy showed that the existence of an $S_2[2,k,n]$ design would imply the existence of an ordinary Steiner system $S(3,2^k,2^n)$. In the case of the parameter set $S_2[2,3,7]$ the corresponding ordinary Steiner system would have the parameters $S(3,8,128)$. Until now, no such Steiner systems have been found. This indicates that it is a really hard problem to find  $S_2[2,3,7]$ designs if they exit. 

In a recent paper \cite{EV12} Etzion and Vardy concentrated on parameters $q=2,t=2,k=3,n=13$. In fact they were able to find a set $\mathcal{B}$ of $14$ orbits on $3$-subspaces of $\F_2^{13}$ with respect to the normalizer of the singer cycle $G\le GL(13,2)$ such that each $2$-subspace is contained in at most one element of $\mathcal{B}$. With this group of automorphisms, a $q$-Steiner system $S_2[2,3,13]$ would consist of $15$ $G$-orbits.

More general, a promising method to find $t$-$(n,k,\lambda;q)$ designs is the well-known Kramer-Mesner method \cite{KM76}. A $t$-$(n,k,\lambda;q)$ design with a prescribed group $G\le GL(n,q)$ of automorphisms is equivalent to a 0-1-solution $x$ of the Diophantine system of equations $A^G x =[\lambda,\ldots,\lambda]^t$ where $A^G$ is the $G$-incidence matrix between the $G$-orbits on $t$-subspaces and $G$-orbits on $k$-subspaces. The entry $a(T,K)$ of $A^G$ whose row corresponds to the $G$-orbit with representative $T$ and whose column corresponds to the $G$-orbit with representative $K$ counts the number of $k$-subspaces in the orbit of $K$ containing the subspace $T$.

Many of all known $t$-$(n,k,\lambda;q)$ designs e.\,g. the design with the largest value $t$ --- a $3-(8,4,12;2)$ design \cite{BKL05} --- were constructed with Kramer-Mesner and used as group of automorphisms the normalizer of the singer cycle, the same group Etzion and Vardy used for their \emph{near}-$q$-Steiner System with $q=2,t=2,k=3,n=13$ \cite{EV12}.

In the following section we give the first construction of a $q$-Steiner system --- a system with parameters $S_2[2,3,13]$ --- which was constructed by Kramer-Mesner using the normalizer of singer cycle as group of automorphisms.

\section{The Construction of $S_2[2,3,13]$}
The normalizer of the Singer cycle $G=\langle F,S\rangle\le GL(13,2)$ is generated by the following two matrices:
\[\footnotesize
F=\left[\begin{array}{ccccccccccccc}
1&0&0&0&0&0&0&0&0&0&0&1&0\\
0&0&0&0&0&0&0&1&0&0&0&1&1\\
0&1&0&0&0&0&0&1&0&0&0&0&0\\
0&0&0&0&0&0&0&0&1&0&0&1&1\\
0&0&1&0&0&0&0&1&1&0&0&1&1\\
0&0&0&0&0&0&0&1&0&1&0&0&0\\
0&0&0&1&0&0&0&0&1&1&0&0&1\\
0&0&0&0&0&0&0&0&1&0&1&0&0\\
0&0&0&0&1&0&0&0&0&1&1&0&0\\
0&0&0&0&0&0&0&0&0&1&0&1&0\\
0&0&0&0&0&1&0&0&0&0&1&1&0\\
0&0&0&0&0&0&0&0&0&0&1&0&1\\
0&0&0&0&0&0&1&0&0&0&0&1&1\\
\end{array}\right]
\]
\[\footnotesize
S=\left[\begin{array}{ccccccccccccc}
0&0&0&0&0&0&0&0&0&0&0&0&1\\
1&0&0&0&0&0&0&0&0&0&0&0&1\\
0&1&0&0&0&0&0&0&0&0&0&0&0\\
0&0&1&0&0&0&0&0&0&0&0&0&1\\
0&0&0&1&0&0&0&0&0&0&0&0&1\\
0&0&0&0&1&0&0&0&0&0&0&0&0\\
0&0&0&0&0&1&0&0&0&0&0&0&0\\
0&0&0&0&0&0&1&0&0&0&0&0&0\\
0&0&0&0&0&0&0&1&0&0&0&0&0\\
0&0&0&0&0&0&0&0&1&0&0&0&0\\
0&0&0&0&0&0&0&0&0&1&0&0&0\\
0&0&0&0&0&0&0&0&0&0&1&0&0\\
0&0&0&0&0&0&0&0&0&0&0&1&0\\
\end{array}\right]
\]
The order of $G$ is $(2^{13}-1)\cdot 13=106483$. The incidence matrix $A^G$ between $2$- and $3$-subspaces has $105$ rows and $25572$ columns whose entries are all $0$ or $1$. Solving the corresponding Diophantine linear system $A^Gx=[1,\ldots,1]^t$ yields at least $26$ solutions which define $q$-Steiner systems $S_2[2,3,13]$. 
The Diophantine linear system has been solved on a standard desktop computer with the dancing links algorithm by D. Knuth \cite{knuth}.
In Table~\ref{tab:1} one solution is given by the $15$ orbit representatives. The union of those $15$ yields a $q$-Steiner system $S_2[2,3,13]$ with $15\cdot106483=1597245$ elements. 

Furthermore, if $A_q(n,d,k)$ denotes the maximum number of codewords of a constant dimension code in $\F_q^n$ with dimension $k$ and minimum subspace distance $d$, we get $A_2(13,4,3)=1597245$.

Since the existence of a $q$-Steiner system $S_2[2,k,n]$ implies the existence of a Steiner system $S(3,2^k,2^n)$ \cite{EV11} we obtain a Steiner system $S(3,8,8192)$ from $S_2[2,3,13]$.

Summarizing these results we can state our main theorem:
\begin{theorem} \
\begin{enumerate}
\item 
$q$-Steiner systems $S_2[2,3,13]$ do exist.

\item
Steiner systems $S(3,8,8192)$ do exist.

\item
The packing bound $A_2(13,4,3)=1597245$ is reached.

\end{enumerate}

\end{theorem}

\begin{table}[!htbp]
\centering
\caption{The $15$ orbit representatives of $S_2[2,3,13]$}\label{tab:1}
\vspace{0.5cm}
{\footnotesize
$\begin{bmatrix}
0&0&0\\
0&0&0\\
0&0&0\\
0&0&0\\
0&0&0\\
1&0&0\\
0&0&0\\
1&0&0\\
1&0&0\\
0&0&0\\
0&0&0\\
0&1&0\\
0&0&1\\
\end{bmatrix}$\qquad
$\begin{bmatrix}
0&0&0\\
0&0&0\\
0&0&0\\
0&0&0\\
0&0&0\\
1&0&0\\
0&0&0\\
0&0&0\\
0&0&0\\
0&0&0\\
0&1&0\\
0&1&0\\
0&0&1\\
\end{bmatrix}$\qquad
$\begin{bmatrix}
0&0&0\\
0&0&0\\
0&0&0\\
0&0&0\\
0&0&0\\
0&0&0\\
1&0&0\\
0&0&0\\
1&0&0\\
0&1&0\\
1&0&0\\
0&0&0\\
0&0&1\\
\end{bmatrix}$\qquad
$\begin{bmatrix}
0&0&0\\
0&0&0\\
0&0&0\\
0&0&0\\
1&0&0\\
1&0&0\\
1&0&0\\
0&0&0\\
1&0&0\\
0&1&0\\
1&1&0\\
1&0&0\\
0&0&1\\
\end{bmatrix}$\qquad
$\begin{bmatrix}
0&0&0\\
0&0&0\\
0&0&0\\
1&0&0\\
0&0&0\\
0&0&0\\
0&0&0\\
0&0&0\\
0&1&0\\
0&0&0\\
0&1&0\\
0&1&0\\
0&0&1\\
\end{bmatrix}$\\[0.5cm]
$\begin{bmatrix}
0&0&0\\
0&0&0\\
1&0&0\\
0&0&0\\
1&0&0\\
0&0&0\\
1&0&0\\
1&0&0\\
0&1&0\\
0&1&0\\
1&1&0\\
1&1&0\\
0&0&1\\
\end{bmatrix}$\qquad
$\begin{bmatrix}
0&0&0\\
0&0&0\\
1&0&0\\
0&0&0\\
0&0&0\\
1&0&0\\
1&0&0\\
0&1&0\\
1&0&0\\
0&0&0\\
1&0&0\\
1&0&0\\
0&0&1\\
\end{bmatrix}$\qquad
$\begin{bmatrix}
1&0&0\\
0&0&0\\
0&0&0\\
1&0&0\\
0&0&0\\
1&0&0\\
1&0&0\\
0&1&0\\
0&0&0\\
1&0&0\\
0&0&0\\
0&1&0\\
0&0&1\\
\end{bmatrix}$\qquad
$\begin{bmatrix}
0&0&0\\
1&0&0\\
1&0&0\\
1&0&0\\
0&0&0\\
0&0&0\\
0&0&0\\
0&1&0\\
0&0&0\\
0&0&0\\
1&1&0\\
0&1&0\\
0&0&1\\
\end{bmatrix}$\qquad
$\begin{bmatrix}
1&0&0\\
0&0&0\\
0&0&0\\
1&0&0\\
1&0&0\\
0&0&0\\
1&0&0\\
0&1&0\\
0&0&0\\
0&1&0\\
1&0&0\\
0&1&0\\
0&0&1\\
\end{bmatrix}$\\[0.5cm]
$\begin{bmatrix}
0&0&0\\
1&0&0\\
0&0&0\\
1&0&0\\
1&0&0\\
1&0&0\\
0&1&0\\
0&0&0\\
1&0&0\\
0&1&0\\
1&0&0\\
0&0&0\\
0&0&1\\
\end{bmatrix}$\qquad
$\begin{bmatrix}
1&0&0\\
0&0&0\\
1&0&0\\
1&0&0\\
1&0&0\\
1&0&0\\
0&1&0\\
0&0&0\\
1&1&0\\
0&0&0\\
1&0&0\\
1&1&0\\
0&0&1\\
\end{bmatrix}$\qquad
$\begin{bmatrix}
1&0&0\\
0&0&0\\
1&0&0\\
1&0&0\\
1&0&0\\
1&0&0\\
0&1&0\\
0&0&0\\
1&1&0\\
0&1&0\\
0&0&0\\
0&1&0\\
0&0&1\\
\end{bmatrix}$\qquad
$\begin{bmatrix}
0&0&0\\
0&0&0\\
0&0&0\\
1&0&0\\
0&0&0\\
1&0&0\\
0&1&0\\
1&1&0\\
1&0&0\\
0&0&0\\
1&0&0\\
0&0&0\\
0&0&1\\
\end{bmatrix}$\qquad
$\begin{bmatrix}
0&0&0\\
1&0&0\\
1&0&0\\
0&0&0\\
0&0&0\\
0&1&0\\
0&0&0\\
0&1&0\\
1&0&0\\
1&0&0\\
0&1&0\\
0&1&0\\
0&0&1\\
\end{bmatrix}$}
\end{table}

\newpage
\section{Acknowledgement}

We would like to thank the attendees of the conference \lq\lq Trends in Coding Theory, Ascona, 2012\rq\rq{} for their fruitful discussions.
\end{document}